\documentclass[12pt]{amsart}

\usepackage{amsmath,xspace,amssymb,mathrsfs}
\usepackage{color}

\input xy
\xyoption{all}
\xyoption{2cell}
\UseAllTwocells
\CompileMatrices

\newcommand{\Spec}{\operatorname{Spec}}
\renewcommand{\phi}{\varphi}

\newcommand{\Ker}{\operatorname{Ker}}

\newcommand{\Min}{\operatorname{Min}}

\newcommand{\Max}{\operatorname{Max}}

\newcommand{\Sp}{\operatorname{Sp}}
\newcommand{\Ann}{\operatorname{Ann}}
\newcommand{\Supp}{\operatorname{Supp}}

\newcommand{\Spp}{\operatorname{Spp}}

\newcommand{\Rad}{\operatorname{Rad}}

\newtheorem{proposition}{Proposition}[section]
\newtheorem{lemma}[proposition]{Lemma} 
\newtheorem{corollary}[proposition]{Corollary}
\newtheorem{theorem}[proposition]{Theorem}

\newtheorem{prop-def}[proposition]{Proposition and definition}

\newtheorem{conjecture}[proposition]{Conjecture}

\theoremstyle{definition}
\newtheorem{definition}[proposition]{Definition}
\newtheorem{example}[proposition]{Example}

\newtheorem{remark}[proposition]{Remark}

\begin{document}

\title[Purely-prime ideals]{On purely-prime ideals with applications}

\author[A. Tarizadeh and M. Aghajani]{Abolfazl Tarizadeh and Mohsen Aghajani}
\address{ Department of Mathematics, Faculty of Basic Sciences,
University of Maragheh \\
P. O. Box 55136-553, Maragheh, Iran.}
\email{ebulfez1978@gmail.com, aghajani14@gmail.com}


\begin{abstract} In this paper, new algebraic and topological results on purely-prime ideals of a commutative ring (pure spectrum) are obtained. Especially, Grothendieck type theorem is obtained which states that there is a canonical correspondence between the idempotents of a ring and the clopens of its pure spectrum. It is also proved that a given ring is a Gelfand ring iff its maximal spectrum equipped with the induced Zariski topology is homeomorphic to its pure spectrum. Then as an application, it is deduced that a ring is zero dimensional iff its prime spectrum and pure spectrum are isomorphic. Dually, it is shown that a given ring is a reduced mp-ring iff its minimal spectrum equipped with the induced flat topology and its pure spectrum are the same.  Finally, the new notion of semi-Noetherian ring is introduced and Cohen type theorem is proved. \\

\end{abstract}

\maketitle

\section{Introduction}

For a given commutative ring $A$, we may assign various spectra
with the canonical maps (see Propositions \ref{Corollary I} and  \ref{Proposition II}): $$\xymatrix{\Spec(A)\ar[r]&\Spp(A)\ar[r]&\Sp(A)}$$
where $\Spec(A)$ is the prime spectrum, $\Spp(A)$ is the pure spectrum and $\Sp(A)$ is the Pierce spectrum whose points are the prime ideals, purely-prime ideals and max-regular ideals, respectively. The prime spectrum is well known, has reach geometric structures and plays a major role in modern algebraic geometry. Two other ones are unknown (or less known) in the literature. It is the purpose of the present paper to study these spectra, specially the pure spectrum, deeply and extensively. Then as a outcome, various interesting results are discovered and some applications are given. \\
Purely-prime ideal notion was introduced and studied in \cite[Chaps. 7 , 8]{Borceux} for general rings (not necessarily commutative),
it is also studied in \cite{Al Ezeh 2}. However, except in these sources, this topic seem to have not been made the subject of special study.  Maybe one of the main reasons that this topic has received less attention over the years is that this natural notion of purely-prime ideal has been unknown (or less known) in the literature. So one of the particular aims of this paper can be considered the introducing of this subject widely to the mathematical community. \\
In this paper we study the purely-prime ideals of a commutative ring and various new and interesting results are obtained. The reason that we focus on the commutative case is that many of the results of this paper do not hold for non-commutative rings. \\
In \S 3 we prove some new and useful (including algebraic and topological) properties of purely-prime ideals of a commutative ring.
One of the main results of this section, Theorem \ref{Theorem III}, establishes a correspondence between the idempotents of a ring and the clopens of its pure spectrum. Pure ideal and purely-prime ideal notions are quite interesting. In sections 4 and 5, we use them in studying and characterizing of Gelfand rings and reduced mp-rings. Especially, it is shown that a ring $A$ is a Gelfand ring iff $\Max(A)\simeq\Spp(A)$, see Theorem \ref{Theorem V}. As an application of this, we obtain that a ring $A$ is zero dimensional iff $\Spec(A)\simeq\Spp(A)$. As another main result, in Theorem \ref{th dual mp-ring}, it is proved that a ring $A$ is a reduced mp-ring iff $\Min(A)=\Spp(A)$ as topological spaces. Theorem \ref{th p.p. new bn} is a further result in this spirit.
In \S 6, we introduce the new notion of semi-Noetherian ring based on the pure ideal notion, and then Cohen type theorem is proved, see Theorem \ref{Theorem Cohen type}. \\ This study also led us to propose a challenging open problem on purely-prime ideals of a commutative ring, see Conjecture \ref{Conjecture I}. \\

\section{Preliminaries}

Most of the following material can be found in \cite[Chaps. 7, 8]{Borceux}. We need them in the sequel and collected here for the convenience of the reader. In this paper, all rings are commutative. Let $I$ be an ideal of a ring $A$. Then $I$ is called a pure ideal if the canonical ring map $A\rightarrow A/I$ is a flat ring map. It is well known that an ideal $I$ of a ring $A$ is a pure ideal iff $\Ann(f)+I=A$ for all $f\in I$, or equivalently, for each $f\in I$ there exists some $g\in I$ such that $f(1-g)=0$, see e.g. \cite[Chap 7, Proposition 2]{Borceux} or \cite[Tag 04PS]{Johan}. If $A$ is either an integral domain or a local ring, then the zero ideal and the whole ring are the only pure ideals of $A$. \\

\begin{lemma}\label{Lemma I} The pure ideals of a ring are stable under taking finite intersections and arbitrary sums. \\
\end{lemma}
{\bf Proof.} \cite[Proposition 1.8]{De Marco}. $\Box$ \\

If $I$ is an ideal of a ring $A$ then we define $\nu(I)$ the sum of all pure ideals of $A$ which are contained in $I$. Such an ideal exists, because the zero ideal is pure. By Lemma \ref{Lemma I}, $\nu(I)$ is the largest pure ideal contained in $I$. It is called the \emph{pure part} of $I$. Clearly $\nu(I\cap J)=\nu(I)\cap\nu(J)$. If $(I_{k})$ is a family of ideals of $A$ then $\sum\limits_{k}\nu(I_{k})\subseteq\nu(\sum\limits_{k}I_{k})$. Later, we shall observe that the equality holds iff $A$ is a Gelfand ring, see Theorem \ref{Theorem V}. Recall that a ring $A$ is said to be a Gelfand ring (or, pm-ring) if each prime ideal of $A$ is contained in a unique maximal ideal of $A$. Dually, a ring $A$ is called a mp-ring if each prime of $A$ contains a unique minimal prime ideal of $A$. It is well known that a ring $A$ is a Gelfand ring if and only if for each maximal ideal $\mathfrak{m}$ of $A$ the canonical ring map $A\rightarrow A_{\mathfrak{m}}$ is surjective, for more information see \cite[Theorem 4.3]{Mohsen-Abolfazl}. \\
If $I$ is an ideal of a ring $A$ then the set $u(I)=\{f\in A: \exists g\in I, f=fg\}$ is an ideal of $A$. It is called the \emph{unit part} of $I$.
Similarly above, we have $u(I\cap J)=u(I)\cap u(J)$ and $\sum\limits_{k}u(I_{k})\subseteq u(\sum\limits_{k}I_{k})$. Moreover, $\nu(I)\subseteq u(I)\subseteq I$. \\

\begin{definition} Every maximal element of the set of proper and pure ideals of a ring $A$ is called a \emph{purely-maximal} ideal of $A$. \\
\end{definition}

By Zorn's Lemma, every non-zero ring has at least a purely-maximal ideal. It follows that every proper and pure ideal of a ring $A$ is contained in a purely-maximal ideal of $A$. \\

\begin{definition} By a \emph{purely-prime} ideal of a ring $A$ we mean a proper and pure ideal $P$ of $A$ such that if there exist pure ideals $I$ and $J$ of $A$ with $IJ\subseteq P$, then either $I\subseteq P$ or $J\subseteq P$. \\
\end{definition}

\begin{lemma}\label{Lemma III} Every purely-maximal ideal is purely-prime. \\
\end{lemma}
{\bf Proof.} See \cite[Chap. 7, Proposition 26]{Borceux}. $\Box$ \\

The set of purely-prime ideals of a ring $A$ is denoted by $\Spp(A)$. If $I$ is a pure ideal of $A$ then we define $U_{I}=\{P\in\Spp(A): I\nsubseteq P\}$. Then clearly $U_{A}=\Spp(A)$ and $U_{I}\cap U_{J}=U_{IJ}$ for pure ideals $I$ and $J$ of $A$ (note that $IJ=I\cap J$ is also a pure ideal of $A$). Thus there exists a (unique) topology over $\Spp(A)$, called the pure topology, such that the base opens are precisely of the form $U_{I}$ where $I$ is a pure ideal of $A$. The set $\Spp(A)$ endowed with the pure topology is called the \emph{pure spectrum} of $A$. Using Lemma \ref{Lemma I}, then it is easy to see that the opens of $\Spp(A)$ are precisely the base opens. In other words, the closed subsets of the pure spectrum $\Spp(A)$ are precisely of the form $V_{p}(I)=\{P\in\Spp(A): I\subseteq P\}$ where $I$ is a pure ideal of $A$. If $P\in\Spp(A)$ then $\overline{\{P\}}=V_{p}(P)$. If $P$ is a purely-prime ideal of $A$ and $e\in A$ is an idempotent, then either $e\in P$ or $1-e\in P$. If $f\in A$ is an idempotent then we shall denote $U_{Af}$ simply by $U_{f}$. \\

\begin{proposition}\label{Proposition I} The pure spectrum $\Spp(A)$ is quasi-compact. \\
\end{proposition}

{\bf Proof.} See \cite[Chap 7, Proposition 34]{Borceux}. $\Box$ \\

Let $I=(f_{i}: i\in S)$ be an ideal of a ring $A$ such that each $f_{i}=f_{i}g_{i}$ for some $g_{i}\in I$. Then $I$ is a pure ideal. \\

\begin{theorem}\label{Theorem nice cvb} If $\phi:A\rightarrow B$ is a morphism of rings then the following hold. \\
$\mathbf{(i)}$ If $I$ is a pure ideal of $A$, then $IB$ is a pure ideal of $B$. \\
$\mathbf{(ii)}$ If $P$ is a purely-prime ideal of $B$, then $\nu\big(\phi^{-1}(P)\big)$ is a purely-prime ideal of $A$. \\
$\mathbf{(iii)}$ The map $\Spp(\phi):\Spp(B)\rightarrow\Spp(A)$ given by $P\rightsquigarrow\nu\big(\phi^{-1}(P)\big)$ is continuous. \\
$\mathbf{(iv)}$ If $\psi:B\rightarrow C$ is a second morphism of rings then $\Spp(\psi\circ\phi)=\Spp(\phi)\circ\Spp(\psi)$. \\
\end{theorem}

{\bf Proof.} For (i) and (ii) see \cite[Chap. 7, Lemmas 60, 62]{Borceux}. For (iii) and (iv) see \cite[Chap 7, Propositions 63, 64]{Borceux}. $\Box$ \\

\begin{lemma}\label{Lemma supp} Let $I$ be an ideal of a ring $A$. Then $\bigcup\limits_{f\in I}D(f)\subseteq\Supp(I)$. The equality holds if and only if $I$ is a pure ideal. \\
\end{lemma}

{\bf Proof.} See \cite[Proposition 1.5]{De Marco}. $\Box$ \\

\begin{lemma}\label{Lemma IV} Let $I$ and $J$ be pure ideals of a ring $A$. Then $U_{I}\subseteq U_{J}$ iff $I\subseteq J$. \\
\end{lemma}

{\bf Proof.} See \cite[Chap. 7, Theorem 32]{Borceux}. $\Box$ \\

If an ideal of a ring $A$ is generated by a set of idempotents of $A$, then it is called a regular ideal of $A$. Each regular ideal is a pure ideal, but the converse does not necessarily hold. Every maximal element of
the set of proper and regular ideals of $A$ is called a max-regular ideal of $A$. The set of max-regular ideals of $A$ is called the \emph{Pierce spectrum} of $A$ and denoted by $\Sp(A)$. It is a compact and totally disconnected topological space whose base opens are precisely of the form $d(f)=\{M\in\Sp(A):f\notin M\}$ where $f\in A$ is an idempotent. For more information see \cite[Lemma 3.18]{Abolfazl}. \\

\section{Pure spectrum versus prime spectrum}

In this section, algebraic and topological properties of the pure spectrum of a commutative ring are studied and various interesting results are obtained. If $\mathfrak{p}$ is a prime ideal of a ring $A$, then the canonical morphism $A\rightarrow A_{\mathfrak{p}}$ is denoted by $\pi_{\mathfrak{p}}$. The following result generalizes \cite[Chap 7, Proposition 27]{Borceux}. \\

\begin{lemma}\label{Lemma II} If $\mathfrak{p}$ is a prime ideal of a ring $A$, then $\nu(\mathfrak{p})$ is a purely-prime ideal of $A$ and $\nu(\mathfrak{p})=\nu(\Ker\pi_{\mathfrak{p}})$. \\
\end{lemma}

{\bf Proof.} Let $I$ and $J$ be two pure ideals of $A$ such that $IJ\subseteq\nu(\mathfrak{p})$. We have then either $I\subseteq\mathfrak{p}$ or $J\subseteq\mathfrak{p}$. Therefore either $I\subseteq\nu(\mathfrak{p})$ or $J\subseteq\nu(\mathfrak{p})$. Therefore $\nu(\mathfrak{p})$ is a purely-prime ideal. Clearly $\Ker\pi_{\mathfrak{p}}\subseteq\mathfrak{p}$ and so $\nu(\Ker\pi_{\mathfrak{p}})\subseteq\nu(\mathfrak{p})$. Conversely,
if $f\in\nu(\mathfrak{p})$ then there exists some $g\in\nu(\mathfrak{p})$ such that $f(1-g)=0$. Clearly $1-g\in A\setminus\mathfrak{p}$. Thus $f\in\Ker\pi_{\mathfrak{p}}$. It follows that $\nu(\mathfrak{p})\subseteq\nu(\Ker\pi_{\mathfrak{p}})$.
$\Box$ \\

\begin{corollary}\label{Corollary IV 2020} If $\mathfrak{p}\subseteq\mathfrak{q}$ are prime ideals of a ring $A$, then $\nu(\mathfrak{p})=\nu(\mathfrak{q})$. \\
\end{corollary}

{\bf Proof.} Clearly $\nu(\mathfrak{p})\subseteq\nu(\mathfrak{q})$. We also have $\Ker\pi_{\mathfrak{q}}\subseteq\Ker\pi_{\mathfrak{p}}$ and so $\nu(\Ker\pi_{\mathfrak{q}})\subseteq\nu(\Ker\pi_{\mathfrak{p}})$. Using Lemma \ref{Lemma II}, then we have $\nu(\mathfrak{q})\subseteq\nu(\mathfrak{p})$. $\Box$ \\

Corollary \ref{Corollary IV 2020}, in particular, tells us that if $M$ is a purely-maximal ideal of a ring $A$, then $M\subseteq\Ker\pi_{\mathfrak{p}}\subseteq\mathfrak{p}$ for some minimal prime ideal $\mathfrak{p}$ of $A$. Note that the converse of Corollary \ref{Corollary IV 2020}, does not necessarily hold. For example, take the ring of integers. \\

\begin{corollary} Let $\mathfrak{p}$ be a prime ideal of a ring $A$ and $I$ a proper ideal of $A$. If $\mathfrak{p}\subseteq I$, then $\nu(\mathfrak{p})=\nu(I)$. \\
\end{corollary}

{\bf Proof.} There exists a maximal ideal $\mathfrak{m}$ of $A$ such that $I\subseteq\mathfrak{m}$. It follows that $\nu(\mathfrak{p})\subseteq\nu(I)\subseteq\nu(\mathfrak{m})$. Then apply Corollary \ref{Corollary IV 2020}. $\Box$ \\

Note that if a prime ideal $\mathfrak{p}$ of a ring $A$ contains an ideal $I$ of $A$. Then the inclusion $\nu(I)\subseteq\nu(\mathfrak{p})$
may be strict. For example, $\mathfrak{p}=\{0,3\}$ is a prime ideal of $A=\mathbb{Z}/6\mathbb{Z}$. We have $\nu(\mathfrak{p})=\mathfrak{p}$, since $3$ is an idempotent. But $\nu(I)=I$ where $I$ is the zero ideal of $A$. \\

If $I$ is a pure ideal of a ring $A$, then its radical $\sqrt{I}$ is not necessarily a pure ideal. For instance, the zero ideal is pure but its radical is not necessarily pure. As a specific example, take $A=\mathbb{Z}/4\mathbb{Z}$ then $\sqrt{0}=\{0,2\}$ is not a pure ideal. \\

\begin{lemma}\label{Lemma Supp 2} Let $I=(f_{i}: i\in S)$ be an ideal of a ring $A$. Then $\Supp(I)=\bigcup\limits_{i}D(f_{i})$ if and only if $I$ is a pure ideal. \\
\end{lemma}

{\bf Proof.} It is proved exactly like Lemma \ref{Lemma supp}. $\Box$ \\

In Lemmas \ref{Corollary V 2020} and \ref{Corollary VI 2020}, we give new and short proofs to \cite[Corollaries 3.4 and 3.5]{Tarizadeh}. \\

\begin{lemma}\label{Corollary V 2020} Let $I$ be an ideal of a ring $A$ such that $\sqrt{I}$ is a pure ideal. Then $I=\sqrt{I}$. \\
\end{lemma}

{\bf Proof.} If $f\in\sqrt{I}$ then there exists some $g\in\sqrt{I}$ such that $f=fg$. Clearly $g^{n}\in I$ for some $n\geqslant1$. We have then $f=fg^{n}\in I$. $\Box$ \\

\begin{lemma}\label{Corollary VI 2020} If $I$ is a pure ideal of a reduced ring $A$, then $I=\sqrt{I}$. \\
\end{lemma}

{\bf Proof.} If $f\in\sqrt{I}$ then there exist a natural number $n\geqslant1$ and some $g\in I$ such that $f^{n}(1-g)=0$. Thus $f(1-g)$ is nilpotent and so $f=fg\in I$. $\Box$ \\

\begin{theorem}\label{Theorem reduced ring} If $I$ is a pure ideal of a reduced ring $A$, then $I=\{f\in A: D(f)\subseteq\Supp(I)\}$. \\
\end{theorem}

{\bf Proof.} Take $f\in A$ such that $D(f)\subseteq\Supp(I)$. By Lemma \ref{Corollary VI 2020}, it suffices to show that $f\in\sqrt{I}$. Suppose there exists a prime ideal $\mathfrak{p}$ of $A$ such that $I\subseteq\mathfrak{p}$ but $f\notin\mathfrak{p}$. Therefore $\mathfrak{p}\in D(f)\subseteq\Supp(I)$ and so there exists some $g\in I$ such that $\Ann(g)\subseteq\mathfrak{p}$. We have $\Ann(g)+I=A$, since $I$ is a pure ideal. But this is a contradiction and we win. $\Box$ \\

Lemma \ref{Lemma II} leads us to a map $\Spec(A)\rightarrow\Spp(A)$ given by $\mathfrak{p}\rightsquigarrow\nu(\mathfrak{p})$. We call it the pure part map and denote by $\nu$. In general, this map is not injective. It seems that it is not also surjective. But it is easy to see that if $M$ is a purely-maximal ideal of $A$ then there exists a maximal ideal $\mathfrak{m}$ of $A$ such that $\nu(\mathfrak{m})=M$. \\

\begin{proposition}\label{Corollary I} The map $\nu:\Spec(A)\rightarrow\Spp(A)$ is continuous. \\
\end{proposition}

{\bf Proof.} Let $I$ be a pure ideal of $A$. If $\mathfrak{p}\in\nu^{-1}(U_{I})$ then $\mathfrak{p}\in\Supp(I)$, because if $I_{\mathfrak{p}}=0$ then $I\subseteq\mathfrak{p}$ and so $I\subseteq\nu(\mathfrak{p})$, which is a contradiction. Therefore by Lemma \ref{Lemma supp},
$\nu^{-1}(U_{I})=\bigcup\limits_{f\in I}D(f)$.  $\Box$ \\

A fundamental result due to Grothendieck states that there is a canonical bijection between the idempotents of a ring and the clopens of its prime spectrum, see \cite[Tag 00EE]{Johan} or \cite[Proposition 3.1]{Tarizadeh 2}. In the following result, we establish the analogue of this theorem for pure spectrum. \\

\begin{theorem}\label{Theorem III}$($Grothendieck type theorem$)$ The map $f\rightsquigarrow U_{f}$ is a bijection from the set of idempotents of a ring $A$ onto the set of clopens of the pure spectrum $\Spp(A)$. \\
\end{theorem}

{\bf Proof.} If $f\in A$ is an idempotent then $Af$ is a pure ideal of $A$. Clearly $U_{f}=V_{p}\big(I)$ where $I=A(1-f)$. Therefore $U_{f}$ is a clopen of $\Spp(A)$ and so the above map is well-defined. Let $f$ and $g$ be two idempotents of $A$ such that $U_{f}=U_{g}$. Then by Lemma \ref{Lemma IV}, $Af=Ag$. Thus $f=ag$ for some $a\in A$. It follows that $f=fg$. Similarly we get that $g=fg$. Thus $f=g$. Then we show that the above map is surjective. Let $F$ be a clopen of the pure spectrum $\Spp(A)$. By Proposition \ref{Corollary I}, $\nu^{-1}(F)$ is a clopen of $\Spec(A)$. Thus there exists an idempotent $e\in A$ such that $\nu^{-1}(F)=D(e)$. We shall prove that $F=U_{e}$. First note that if $P$ is a purely prime ideal of $A$ then there exists a purely maximal ideal $M$ of $A$ such that $P\subseteq M$. We also have $M=\nu(\mathfrak{m})\subseteq\mathfrak{m}$ for some $\mathfrak{m}\in\Max(A)$. Now if $P\in F$ then $M\in V_{p}(P)=\overline{\{P\}}\subseteq F$. Thus $\mathfrak{m}\in D(e)$ and so $e\notin P$. Hence, $F\subseteq U_{e}$. Conversely, if $P\in U_{e}$ then $1-e\in P$. Thus $\mathfrak{m}\in D(e)$ and so $M\in F$. But $F$ is an open set, so $F=U_{I}$ for some (pure) ideal $I$. Thus $I\nsubseteq M$.  This yields that $I\nsubseteq P$. Hence, $P\in F$. $\Box$ \\

\begin{corollary}\label{Coro conn 1089} Let $A$ be a ring. Then the pure spectrum $\Spp(A)$ is connected iff $A$ has no nontrivial idempotents. \\
\end{corollary}

{\bf Proof.} It follows from Theorem \ref{Theorem III}. $\Box$ \\

\begin{proposition}\label{Proposition II} The map
$\lambda:\Spp(A)\rightarrow\Sp(A)$ given by $P\rightsquigarrow(f\in P: f=f^{2})$ is continuous and surjective. \\
\end{proposition}

{\bf Proof.} If $P$ is a purely-prime ideal of $A$ then it is easy to see that $\lambda(P)$ is a max-regular ideal of $A$. Thus $\lambda$ is well-defined.
If $f\in A$ is an idempotent then clearly $\lambda^{-1}\big(d(f)\big)=U_{f}$. Hence $\lambda$ is continuous. Finally, let $M$ be a max-regular ideal of $A$. There exists a maximal ideal $\mathfrak{m}$ of $A$ such that $M\subseteq\mathfrak{m}$. Clearly $M\subseteq\nu(\mathfrak{m})$ and so $M\subseteq\lambda\big(\nu(\mathfrak{m})\big)$. It follows that $M=\lambda\big(\nu(\mathfrak{m})\big)$. $\Box$ \\

\begin{lemma}\label{Lemma v 2020 bn} If $I$ is a pure ideal of a ring $A$, then the following hold. \\
$\mathbf{(i)}$ The pure ideals of $A/I$ are precisely of the form $J/I$ where $J$ is a pure ideal of $A$ such that $I\subseteq J$. \\
$\mathbf{(ii)}$ $\Spp(A/I)=\{P/I: P\in V_{p}(I)\}$. \\
\end{lemma}

{\bf Proof.} $\mathbf{(i)}$ If $J$ is a pure ideal of $A$ such that $I\subseteq J$ then by Theorem \ref{Theorem nice cvb}(i), $J/I$ is a pure ideal. Conversely, if $J/I$ is a pure ideal of $A/I$ then we show that $J$ is a pure ideal. If $f\in J$ then there is some $g\in J$ such that $f(1-g)\in I$. Thus $f(1-g)(1-h)=0$ for some $h\in I$. So $f(1-x)=0$ where $x:=g+h-gh\in J$. \\
$\mathbf{(ii)}$ It follows from (i) and the definition of purely-prime ideal.
$\Box$ \\

\begin{corollary}\label{Lemma VI} If $I$ is a pure ideal of a ring $A$, then the canonical ring map $\pi: A\rightarrow A/I$ induces a homeomorphism from the pure spectrum $\Spp(A/I)$ onto $V_{p}(I)$. \\
\end{corollary}

{\bf Proof.} Using Lemma \ref{Lemma v 2020 bn} and Theorem \ref{Theorem nice cvb}(iii), then the map $\Spp(\pi):\Spp(A/I)\rightarrow\Spp(A)$ given by $P/I\rightsquigarrow\nu(P)=P$ is continuous, injective and its image is $V_{p}(I)$. It is also a closed map, because if $F$ is a closed subset of $\Spp(A/I)$ then by Lemma \ref{Lemma v 2020 bn}, $F=V_{p}(J/I)$ where $J$ is a pure ideal, and we have $\Spp(\pi)(F)=V_{p}(J)$ which is a closed subset of $\Spp(A)$. $\Box$ \\

\begin{example} It is important to notice that if $I$ is not a pure ideal then Corollary \ref{Lemma VI} does not hold. In fact, $\Spp(A/I)$ is not necessarily homeomorphic to $V_{p}\big(\nu(I)\big)$.
As an example, let $m\geq2$ be an integer with the prime factorization $m=p^{c_{1}}_{1}...p^{c_{k}}_{k}$ where the $p_{i}$ are distinct prime numbers and $c_{i}\geq1$ for all $i$. We have $\nu(p_{i}\mathbb{Z}/m\mathbb{Z})=
p^{c_{i}}_{i}\mathbb{Z}/m\mathbb{Z}$
and $\Spp(\mathbb{Z}/m\mathbb{Z})=
\{p^{c_{1}}_{1}\mathbb{Z}/m\mathbb{Z},...,
p^{c_{k}}_{k}\mathbb{Z}/m\mathbb{Z}\}$ but $V_{p}\big(\nu(m\mathbb{Z})\big)=V_{p}(0)=
\Spp(\mathbb{Z})=\{0\}$. \\
\end{example}

\begin{theorem}\label{Theorem 98 connected comp.} The connected components of the pure spectrum $\Spp(A)$ are precisely of the form $V_{p}(M)$ where $M$ is a max-regular ideal of $A$. \\
\end{theorem}

{\bf Proof.} If $C$ is a connected component of $\Spp(A)$, then there exists a max-regular ideal $M$ of $A$ such that $\lambda(C)=\{M\}$ because $\Sp(A)$ is totally disconnected and for $\lambda$ see Proposition \ref{Proposition II}. It follows that $C\subseteq\lambda^{-1}(\{M\})=V_{p}(M)$. By Corollary \ref{Lemma VI}, $V_{p}(M)$ is homeomorphic to $\Spp(A/M)$. But $A/M$ has no nontrivial idempotents, see \cite[Lemma 3.19]{Abolfazl}. Thus by Corollary \ref{Coro conn 1089}, $V_{p}(M)$
a connected subset of $\Spp(A)$. Therefore $C=V_{p}(M)$. Conversely, if $M$ is a max-regular ideal of $A$ then, in the above, we observed that $V_{p}(M)$ is a connected subset of $\Spp(A)$. So it is contained in a connected component $C'$ of $\Spp(A)$. By what we have just proven, there exists a max-regular ideal $M'$ of $A$ such that $C'=V_{p}(M')$. By Lemma \ref{Lemma IV}, $M'\subseteq M$. It follows that $M'=M$. $\Box$ \\

By $\pi_{0}(X)$ we mean the space of connected components of a space $X$. \\

\begin{corollary} $\pi_{0}\big(\Spp(A)\big)$ is canonically homeomorphic to $\Sp(A)$. \\
\end{corollary}

{\bf Proof.} It follows from Theorem \ref{Theorem 98 connected comp.}. $\Box$ \\

\begin{remark}\label{Theorem purely-minimal} By a \emph{purely-minimal} ideal of a ring $A$ we mean a purely-prime ideal $P$ of $A$ such that if there exists a purely-prime ideal $P'$ of $A$ with $P'\subseteq P$, then $P'=P$. If $P$ is a purely-prime ideal of $A$, then there exists a purely-minimal ideal of $A$ contained in $P$, because
if $\mathcal{S}$ is the set of all purely-prime ideals of $A$ which are contained in $P$, then clearly it is a non-empty set, and if $\mathcal{C}$ is a chain in $\mathcal{S}$, then it is easy to see that $\nu(\bigcap\limits_{P'\in\mathcal{C}}P')$ is a purely-prime ideal of $A$, therefore by the Zorn's lemma, $\mathcal{S}$ has at least a minimal element. \\
\end{remark}

\section{Characterizations of Gelfand rings}

In this section new characterizations of Gelfand rings based on ``pure part'' and ``unit part'' notions are given. The following result characterizes the purely-maximal ideals of a Gelfand ring. \\

\begin{theorem}\label{Proposition IV} The purely-maximal ideals of a Gelfand ring $A$ are precisely of the form $\Ker\pi_{\mathfrak{m}}$ where $\mathfrak{m}$ is a maximal ideal of $A$. \\
\end{theorem}

{\bf Proof.} If $\mathfrak{m}$ is a maximal ideal of $A$, then by \cite[Theorem 4.3(vi)]{Mohsen-Abolfazl}, $A/\Ker\pi_{\mathfrak{m}}$ is canonically isomorphic to $A_{\mathfrak{m}}$. Thus $A/\Ker\pi_{\mathfrak{m}}$ is a flat $A-$module and so $\Ker\pi_{\mathfrak{m}}$ is a pure ideal. Therefore there exists a purely-maximal ideal $M$ of $A$ such that $\Ker\pi_{\mathfrak{m}}\subseteq M$. There exists a maximal ideal $\mathfrak{m}'$ of $A$ such that $M=\nu(\mathfrak{m}')\subseteq\Ker\pi_{\mathfrak{m}'}\subseteq
\mathfrak{m}'$. If $\mathfrak{m}\neq\mathfrak{m}'$ then by \cite[Theorem 4.3(ix)]{Mohsen-Abolfazl}, we have
$\Ker\pi_{\mathfrak{m}}+\Ker\pi_{\mathfrak{m}'}=A$. But this is a contradiction. Therefore $\Ker\pi_{\mathfrak{m}}=M$. $\Box$ \\

In Theorem \ref{Theorem V}, it is shown that every purely-prime ideal of a Gelfand ring is purely-maximal. In order to prove Theorem \ref{Theorem V}, the whole strength of Theorem \ref{Theorem IX} and Corollary \ref{Theorem VI} are used. \\

\begin{theorem}\label{Theorem IX} Let $A$ be a Gelfand ring. If the unit part of a maximal ideal $\mathfrak{m}$ of $A$ is contained in a proper ideal $I$ of $A$, then $I\subseteq\mathfrak{m}$. \\
\end{theorem}

{\bf Proof.} See \cite[Chap 8, Lemma 13]{Borceux}. $\Box$ \\

\begin{corollary}\label{Theorem VI} Let $A$ be a Gelfand ring and $I$ be an ideal of $A$. If $u(I)\subseteq\mathfrak{m}$ for some maximal ideal $\mathfrak{m}$ of $A$, then $I\subseteq\mathfrak{m}$. In particular, $\nu(I)=u(I)$. \\
\end{corollary}

{\bf Proof.} See \cite[Chap 8, Propositions 29, 30]{Borceux}. $\Box$ \\

\begin{lemma}\label{Lemma VIII} If $I$ and $J$ are ideals of a ring $A$, then $\nu(IJ)=\nu(I)\nu(J)$. \\
\end{lemma}

{\bf Proof.} We have $IJ\subseteq I\cap J$. It follows that $\nu(IJ)\subseteq\nu(I\cap J)=\nu(I)\cap\nu(J)=\nu(I)\nu(J)\subseteq IJ$. But $\nu(IJ)$ is the largest pure ideal contained in $IJ$, therefore $\nu(IJ)=\nu(I)\nu(J)$. $\Box$ \\

\begin{corollary} If $I$ and $J$ are ideals of a ring $A$, then $u(I)u(J)\subseteq u(IJ)$. If moreover $A$ is a Gelfand ring, then the equality holds.  \\
\end{corollary}

{\bf Proof.} To see the first inclusion it suffices to show that if $f\in u(I)$ and $g\in u(J)$ then $fg\in u(IJ)$. There exist $f'\in I$ and $g'\in J$ such that $f=ff'$ and $g=gg'$ and so $fg=fg(f'g')$. Thus $fg\in u(IJ)$. If $A$ is a Gelfand ring then by Corollary \ref{Theorem VI} and Lemma \ref{Lemma VIII}, $u(IJ)=\nu(IJ)=\nu(I)\nu(J)=u(I)u(J)$. $\Box$ \\

If $I$ is an ideal of a ring $A$ then the intersection of all maximal ideals
of $A$ containing $I$ is denoted by $\Rad(I)$. In Theorem \ref{Theorem V}, we have improved \cite[Chap 8, Theorem 31]{Borceux} and \cite[Theorem 3.7]{Borceux 2} by adding (v)-(viii) as new equivalents. This theorem characterizes Gelfand rings in terms of the pure parts of ideals. Then in Theorem \ref{Theorem VII}, further characterizations of Gelfand rings are given where the unit parts of ideals are involved. \\

\begin{theorem} $($pure part characterization$)$\label{Theorem V} For a ring $A$ the following statements are equivalent. \\
$\mathbf{(i)}$ $A$ is a Gelfand ring. \\
$\mathbf{(ii)}$ If $I$ and $J$ are coprime ideals of $A$ then so are $\nu(I)$ and $\nu(J)$. \\
$\mathbf{(iii)}$ If $(I_{k})$ is a family of ideals of $A$ then $\nu(\sum\limits_{k}I_{k})=\sum\limits_{k}\nu(I_{k})$. \\
$\mathbf{(iv)}$ If $I$ is an ideal of $A$ then $\Rad(I)=\Rad\big(\nu(I)\big)$. \\
$\mathbf{(v)}$ If $\mathfrak{m}$ and $\mathfrak{m}'$ are distinct maximal ideals of $A$ then $\nu(\mathfrak{m})+\nu(\mathfrak{m}')=A$. \\
$\mathbf{(vi)}$ If the pure part of an ideal $I$ of $A$ is contained in a maximal ideal $\mathfrak{m}$ of $A$, then $I\subseteq\mathfrak{m}$. \\
$\mathbf{(vii)}$ If $I$ is an ideal of $A$, then $\Max(A)\cap V(I)=\Max(A)\cap V\big(\nu(I)\big)$. \\
$\mathbf{(viii)}$ The map $\Max(A)\rightarrow\Spp(A)$ given by $\mathfrak{m}\rightsquigarrow\nu(\mathfrak{m})$ is an isomorphism. \\
\end{theorem}

{\bf Proof.} For $\mathbf{(i)}\Leftrightarrow\mathbf{(ii)}\Leftrightarrow\mathbf{(iii)}$ see \cite[Chap 8, Theorem 31]{Borceux}. For $\mathbf{(i)}\Leftrightarrow\mathbf{(iv)}$ see \cite[Theorem 3.7]{Borceux 2}. \\
$\mathbf{(i)}\Rightarrow\mathbf{(v)}:$ If $\mathfrak{m}$ is a maximal ideal of $A$ then $\Ker\pi_{\mathfrak{m}}$ is a pure ideal and so by Lemma \ref{Lemma II}, $\nu(\mathfrak{m})=\Ker\pi_{\mathfrak{m}}$. For distinct maximal ideals $\mathfrak{m}$ and $\mathfrak{m}'$ we have $\Ker\pi_{\mathfrak{m}}+\Ker\pi_{\mathfrak{m}'}=A$. \\
$\mathbf{(v)}\Rightarrow\mathbf{(i)}:$ If $\mathfrak{m}$ and $\mathfrak{m}'$ are distinct maximal ideals of $A$ then there exist $x\in\nu(\mathfrak{m})$ and $y\in\nu(\mathfrak{m}')$ such that $x+y=1$. There are also $f\in\Ann(x)$ and $g\in\mathfrak{m}$ such that $f+g=1$. Clearly $f\in A\setminus\mathfrak{m}$ and $x\in A\setminus\mathfrak{m}'$. $\mathbf{(i)}\Rightarrow\mathbf{(vi)}:$ See Corollary \ref{Theorem VI}. The implications
$\mathbf{(vi)}\Rightarrow\mathbf{(vii)}\Rightarrow\mathbf{(iv)}$ are easy. $\mathbf{(i)}\Rightarrow\mathbf{(viii)}:$ The above map by Proposition \ref{Corollary I} is continuous and by (v) is injective. To see its surjectivity it suffices to show that every purely-prime ideal of $A$ is purely-maximal, see also \cite[Chap 8, Proposition 37]{Borceux}. If $P$ is a purely-prime ideal of $A$ then there exists a purely-maximal ideal $M$ of $A$ such that $P\subseteq M$. If the inclusion is strict then by (iii), there exists some $f\in M$ such that $\nu(Af)\nsubseteq P$. But we have $\nu(Af)\nu\big(\Ann(f)\big)\subseteq(Af)\Ann(f)=0$. It follows that $\nu\big(\Ann(f)\big)\subseteq P$. We also have $\Ann(f)+M=A$. Then by (ii), we get that $\nu\big(\Ann(f)\big)+M=A$. But this is a contradiction. Thus $P$ is a purely-maximal ideal of $A$. Then we show that the pure spectrum $\Spp(A)$ is Hausdorff, see also \cite[Chap 8, Theorem 39]{Borceux}. If $M$ and $M'$ are distinct purely-maximal ideals of $A$ then there exist distinct maximal ideals $\mathfrak{m}$ and $\mathfrak{m}'$ of $A$ such that $M=\nu(\mathfrak{m})$ and $M'=\nu(\mathfrak{m}')$. By the hypothesis, there exist $f\in A\setminus\mathfrak{m}$ and $g\in A\setminus\mathfrak{m}'$ such that $fg=0$. Take $I=\nu(Af)$. Then $M\in U_{I}$, because $\mathfrak{m}+Af=A$ and so $M+I=A$. Similarly, $M'\in U_{J}$ where $J=\nu(Ag)$. But $IJ\subseteq(Af)(Ag)=0$ and so $U_{I}\cap U_{J}=\emptyset$. Therefore $\Spp(A)$ is a Hausdorff space. Moreover, for any ring $A$, then $\Max(A)$ is quasi-compact in the Zariski topology. Hence, the above map is a closed map. $\mathbf{(viii)}\Rightarrow\mathbf{(i)}:$ By the hypothesis and Proposition \ref{Corollary I}, $\Max(A)$ is the Zariski retraction of $\Spec(A)$. Thus by \cite[Theorem 4.3(ii)]{Mohsen-Abolfazl}, $A$ is a Gelfand ring. $\Box$ \\

In particular, we obtain the following nontrivial results. \\

\begin{corollary} If $I$ is a pure ideal of a Gelfand ring $A$, then $I=\sum\limits_{f\in I}\nu(Af)$. \\
\end{corollary}

{\bf Proof.} We have $I=\sum\limits_{f\in I}Af$. Thus by Theorem \ref{Theorem V}(iii), $I=\nu(I)=\sum\limits_{f\in I}\nu(Af)$. $\Box$ \\

\begin{corollary} Let $A$ be a ring. Then $A$ is zero dimensional iff the pure part map $\nu:\Spec(A)\rightarrow\Spp(A)$ is an isomorphism. \\
\end{corollary}

{\bf Proof.} Every zero dimensional ring is a Gelfand ring, so the implication
``$\Rightarrow$'' is deduced from Theorem \ref{Theorem V}(viii).
To see the converse, it will be enough to show that every maximal ideal of $A$ is the radical of a pure ideal. In fact, we shall prove that $\mathfrak{m}=\sqrt{\nu(\mathfrak{m})}$ for all $\mathfrak{m}\in\Max(A)$. Let $\mathfrak{p}$ be a prime ideal of $A$ such that $\nu(\mathfrak{m})\subseteq\mathfrak{p}$. It follows that $\nu(\mathfrak{m})\subseteq\nu(\mathfrak{p})$. Thus $\nu(\mathfrak{p})\in V_{p}\big(\nu(\mathfrak{m})\big)=
\overline{\{\nu(\mathfrak{m})\}}$. This yields that $\mathfrak{p}\in\overline{\{\mathfrak{m}\}}=\{\mathfrak{m}\}$, since $\nu^{-1}$ is continuous. So $\mathfrak{p}=\mathfrak{m}$. $\Box$ \\

\begin{remark} Note that if $A$ is a zero dimensional ring, then by \cite[Corollary 3.17]{Tarizadeh-Aghajani 2}, $\nu(\mathfrak{p})=\Ker\pi_{\mathfrak{p}}$ for all $\mathfrak{p}\in\Spec(A)$. \\
\end{remark}

\begin{theorem} $($unit part characterization$)$\label{Theorem VII} For a ring $A$ the following statements are equivalent. \\
$\mathbf{(i)}$ $A$ is a Gelfand ring. \\
$\mathbf{(ii)}$ If the unit part of an ideal $I$ of $A$ is contained in a maximal ideal $\mathfrak{m}$ of $A$, then $I\subseteq\mathfrak{m}$. \\
$\mathbf{(iii)}$ If $\mathfrak{m}$ and $\mathfrak{m}'$ are distinct maximal ideals of $A$, then $u(\mathfrak{m})+u(\mathfrak{m}')=A$. \\
$\mathbf{(iv)}$ If $I$ and $J$ are coprime ideals of $A$ then so are $u(I)$ and $u(J)$. \\
$\mathbf{(v)}$ If $(I_{k})$ is a family of ideals of $A$, then $u(\sum\limits_{k}I_{k})=\sum\limits_{k}u(I_{k})$. \\
$\mathbf{(vi)}$ If $I$ is an ideal of $A$, then $\Max(A)\cap V(I)=\Max(A)\cap V\big(u(I)\big)$. \\
$\mathbf{(vii)}$ If $I$ is an ideal of $A$, then $\Rad(I)=\Rad\big(u(I)\big)$. \\
$\mathbf{(viii)}$ The map $\Max(A)\rightarrow\Spp(A)$ given by $\mathfrak{m}\rightsquigarrow u(\mathfrak{m})$ is an isomorphism. \\
\end{theorem}

{\bf Proof.} $\mathbf{(i)}\Rightarrow\mathbf{(ii)}:$ See Corollary \ref{Theorem VI}. $\mathbf{(ii)}\Rightarrow\mathbf{(iii)}:$ Easy. \\
$\mathbf{(iii)}\Rightarrow\mathbf{(i)}:$ If $\mathfrak{m}$ and $\mathfrak{m}'$ are distinct maximal ideals of $A$ then by the hypothesis, there exist $a\in u(\mathfrak{m})$ and $b\in u(\mathfrak{m}')$ such that $a+b=1$. There exist also $f\in\mathfrak{m}$ and $g\in\mathfrak{m}'$ such that $a(1-f)=0$ and $b(1-g)=0$. Then clearly $1-f\in A\setminus\mathfrak{m}$, $1-g\in A\setminus\mathfrak{m}'$ and $(1-f)(1-g)=0$. Thus by \cite[Theorem 4.3(v)]{Mohsen-Abolfazl}, $A$ is a Gelfand ring. $\mathbf{(ii)}\Rightarrow\mathbf{(iv)}:$ Easy. 
$\mathbf{(iv)}\Rightarrow\mathbf{(v)}:$ Clearly $\sum\limits_{k}u(I_{k})\subseteq u(\sum\limits_{k}I_{k})$. Conversely, if $f\in u(\sum\limits_{k}I_{k})$ then there exists some $g\in\sum\limits_{k}I_{k}$ such that $f=fg$. We have $g\in\sum\limits_{k=1}^{n}I_{k}$. It follows that $\sum\limits_{k=1}^{n}I_{k}+\Ann(f)=A$. Then by the iteration using of the hypothesis, we get that $\sum\limits_{k=1}^{n}u(I_{k})+\Ann(f)=A$. Hence there exist $h\in\sum\limits_{k=1}^{n}u(I_{k})$ and $h'\in\Ann(f)$ such that $h+h'=1$. This yields that $f=fh\in\sum\limits_{k=1}^{n}u(I_{k})\subseteq\sum\limits_{k}u(I_{k})$. The implications $\mathbf{(v)}\Rightarrow\mathbf{(iii)}$ and $\mathbf{(ii)}\Rightarrow\mathbf{(vi)}\Rightarrow\mathbf{(vii)}
\Rightarrow\mathbf{(ii)}$ are easy. $\mathbf{(i)}\Rightarrow\mathbf{(viii)}:$ If $A$ is a Gelfand ring then the unit and pure parts are the same, then apply
Theorem \ref{Theorem V} (viii). $\mathbf{(viii)}\Rightarrow\mathbf{(i)}:$ If $\mathfrak{m}$ is a maximal ideal of $A$ then by the hypothesis, $u(\mathfrak{m})$ is a pure ideal and so $\nu(\mathfrak{m})=u(\mathfrak{m})$. Then apply Theorem \ref{Theorem V} (viii). $\Box$ \\

The following result generalizes \cite[Chap. 8, Proposition 17]{Borceux} to any ring. \\
\begin{proposition} If $I$ is an ideal of a ring $A$, then $u(I)=\bigcap\limits_{\mathfrak{m}\in\Max(A)\cap V(I)}u(\mathfrak{m})$. \\
\end{proposition}
{\bf Proof.} If $f\in\bigcap\limits_{\mathfrak{m}\in\Max(A)\cap V(I)}u(\mathfrak{m})$ then $\Ann(f)+I=A$. It follows that $f\in u(I)$. $\Box$ \\

\section{mp-rings}

\begin{lemma}\label{Lemma IX 2020} If $I$ and $J$ are pure ideals of a ring $A$ such that $\sqrt{I}=\sqrt{J}$, then $I=J$. \\
\end{lemma}

{\bf Proof.} If $f\in I$ then $f=fg$ for some $g\in I$. Clearly $g^{n}\in J$ for some $n\geqslant1$. Thus $f=fg^{n}\in J$. Hence, $I\subseteq J$. Similarly we get that $J\subseteq I$. $\Box$ \\

Although ``pure part'' and ``unit part'' are quite efficient tools in characterizing Gelfand rings, but similar statements do not the right criteria for characterizing mp-rings. For instance, if $p$ and $q$ are distinct prime numbers then $p\mathbb{Z}+q\mathbb{Z}=\mathbb{Z}$ but $\nu(p\mathbb{Z})+\nu(q\mathbb{Z})=0$. In spite of this, we have the following result. \\

\begin{theorem}\label{Theorem I 2020} Let $A$ be a ring. Then $A$ is a mp-ring if and only if $\mathfrak{p}=\sqrt{\nu(\mathfrak{p})}$ for all $\mathfrak{p}\in\Min(A)$. \\
\end{theorem}

{\bf Proof.} If $A$ is a mp-ring, then by \cite[Theorem 6.2(x)]{Mohsen-Abolfazl}, there exists a pure ideal $I$ of $A$ such that $\mathfrak{p}=\sqrt{I}$. Thus $I\subseteq\nu(\mathfrak{p})$ and so
$\mathfrak{p}=\sqrt{I}\subseteq\sqrt{\nu(\mathfrak{p})}
\subseteq\mathfrak{p}$. Therefore
$\sqrt{I}=\sqrt{\nu(\mathfrak{p})}$. Thus by Lemma \ref{Lemma IX 2020}, $I=\nu(\mathfrak{p})$. For the reverse implication see
\cite[Theorem 6.2(x)]{Mohsen-Abolfazl}. $\Box$ \\

The following result characterizes the purely-maximal ideals of a  mp-ring. \\

\begin{theorem}\label{Theorem II 2020} The purely-maximal ideals of a mp-ring $A$ are precisely of the form $\nu(\mathfrak{p})$ where $\mathfrak{p}$ is a minimal prime of $A$. \\
\end{theorem}

{\bf Proof.} If $M$ is a purely-maximal ideal of $A$, then there exists a maximal ideal $\mathfrak{m}$ of $A$ such that $M=\nu(\mathfrak{m})$. There exists a minimal prime $\mathfrak{q}$ of $A$ such that $\mathfrak{q}\subseteq\mathfrak{m}$. Then by Corollary \ref{Corollary IV 2020}, $M=\nu(\mathfrak{q})$. Conversely, if $\mathfrak{p}$ is a minimal prime of $A$, then there exists a purely-maximal ideal $M$ of $A$ such that $\nu(\mathfrak{p})\subseteq M$. We observed that $M=\nu(\mathfrak{q})$ for some $\mathfrak{q}\in\Min(A)$. So $\nu(\mathfrak{p})\subseteq\mathfrak{q}$. Thus by Theorem \ref{Theorem I 2020}, $\mathfrak{p}=\mathfrak{q}$. $\Box$ \\

\begin{corollary}\label{Proposition III} The purely-maximal ideals of a reduced mp-ring $A$ are precisely the minimal primes of $A$. \\
\end{corollary}

{\bf Proof.} If $\mathfrak{p}$ is a minimal prime of $A$, then by Theorem \ref{Theorem I 2020}, $\mathfrak{p}=\sqrt{\nu(\mathfrak{p})}$. Thus by Lemma \ref{Corollary VI 2020}, $\mathfrak{p}=\nu(\mathfrak{p})$. Then apply Theorem \ref{Theorem II 2020}. $\Box$ \\

Then we give a new proof to the following result. \\

\begin{theorem}\cite[Theorem 3.5]{Al Ezeh 2}\label{Theorem Hasan} Every purely-prime ideal of a reduced mp-ring $A$ is purely-maximal. \\
\end{theorem}

{\bf Proof.} Let $P$ be a purely-prime ideal of $A$. There exists a purely-maximal ideal $M$ of $A$ such that $P\subseteq M$. If $f\in M$ then clearly $J_{f}\cap\Ann(f)$ is contained in the Jacobson radical of $A$ where $J_{f}:=\bigcap\limits_{\mathfrak{m}\in\Max(A)\cap V(f)}\Ker\pi_{\mathfrak{m}}$. By \cite[Theorem 6.4(iv) and Corollary 7.4]{Mohsen-Abolfazl}, $\Ann(f)$ and $J_{f}$ are pure ideals. Therefore $J_{f}\cap\Ann(f)=0$. It follows that $J_{f}\subseteq P$ because $\Ann(f)+M=A$. There exists some $h\in M$ such that $f(1-h)=0$. This yields that $f\in J_{h}$. But, as we observed in the above, $J_{h}\subseteq P$. Therefore $P=M$. $\Box$ \\

Let $A$ be a ring. There exists a (unique) topology over $\Spec(A)$, called the flat topology, such that the collection of $V(I)$, with $I$ is a finitely generated ideals of $A$, forms a base for its opens. For more information see \cite{Abolfazl}. Then we have the following result. \\

\begin{theorem}\label{th dual mp-ring} Let $A$ be a ring and consider the induced flat topology over $\Min(A)$. Then the following statements are equivalent. \\
$\mathbf{(i)}$ $A$ is a reduced mp-ring. \\
$\mathbf{(ii)}$ $\Min(A)=\Spp(A)$ as topological spaces. \\
$\mathbf{(iii)}$ $\Min(A)=\Spp(A)$ as sets. \\
\end{theorem}

{\bf Proof.} $\mathbf{(i)}\Rightarrow\mathbf{(ii)}:$ If $A$ is a reduced mp-ring then by Corollary \ref{Proposition III} and Theorem \ref{Theorem Hasan}, $\Min(A)=\Spp(A)$ as sets. Then we show that they are equal as topological spaces. If $f\in A$ then by \cite[Theorem 6.4(iv)]{Mohsen-Abolfazl}, $I=\Ann(f)$ is a pure ideal. We have $\Min(A)\cap V(f)=U_{I}$. But for any ring $A$, the collection of $\Min(A)\cap V(f)$ with $f\in A$ forms a subbase for the opens of the induced flat topology over $\Min(A)$. Thus the pure topology over $\Spp(A)$ is finer than the induced flat topology. The induced flat topology over $\Min(A)$ is Hausdorff, because if $\mathfrak{p}$ and $\mathfrak{q}$ are distinct minimal prime ideals of $A$ then $\mathfrak{p}+\mathfrak{q}=A$ and so there are $f\in\mathfrak{p}$ and $g\in\mathfrak{q}$ such that $f+g=1$, this yields that $V(f)\cap V(g)=\emptyset$. By Proposition \ref{Proposition I}, the pure topology over $\Spp(A)$ is quasi-compact. Therefore these two topologies are the same, because it is well known that if $\mathscr{T}\subseteq\mathscr{T}'$ are two topologies over a set $X$ such that $\mathscr{T}$ is Hausdorff and $\mathscr{T}'$ is quasi-compact, then $\mathscr{T}=\mathscr{T}'$. $\mathbf{(ii)}\Rightarrow\mathbf{(iii)}:$ There is nothing to prove.
$\mathbf{(iii)}\Rightarrow\mathbf{(i)}:$ Each minimal prime of $A$ is a pure ideal. If  $f\in A$ is a nilpotent then $\Ann(f)=A$, because if it is a proper ideal then it is contained in a maximal ideal $\mathfrak{m}$ of $A$, there exists a minimal prime $\mathfrak{p}$ of $A$ such that $\mathfrak{p}\subseteq\mathfrak{m}$, but $f\in\mathfrak{p}$ and so $\Ann(f)+\mathfrak{p}=A$ which is a contradiction. Thus $f=0$. Hence, $A$ is a reduced ring. If $\mathfrak{p}$ and $\mathfrak{q}$ are minimal primes of $A$ contained in a maximal ideal of $A$, then for each $f\in\mathfrak{p}$ there is some $g\in\mathfrak{p}$ such that $f(1-g)=0$, but $1-g\notin\mathfrak{q}$ and so $f\in\mathfrak{q}$, thus $\mathfrak{p}=\mathfrak{q}$. Therefore $A$ is a mp-ring. $\Box$ \\

Recall that a ring $A$ is called a p.p. ring if every principal ideal of $A$ is a projective $A-$module, or equivalently, $\Ann(f)$ is generated by an idempotent element for all $f\in A$. \\

\begin{theorem}\label{th p.p. new bn} Let $A$ be a reduced mp-ring. Then $A$ is a p.p. ring iff the pure topology and the induced Zariski topology over the set $\Min(A)=\Spp(A)$ are the same. \\
\end{theorem}

{\bf Proof.} Let $A$ be a p.p. ring. Clearly the induced Zariski topology over $\Min(A)$ is finer than the pure topology. To see the reverse inclusion, take $\mathfrak{p}\in\Min(A)\cap D(f)$ where $f\in A$. There is an idempotent $e\in A$ such that $\Ann(f)=Ae$. We have $\mathfrak{p}\in U_{1-e}\subseteq\Min(A)\cap D(f)$, because if $\mathfrak{q}\in U_{1-e}$ then $e\in\mathfrak{q}$ and so $\Ann(f)\subseteq\mathfrak{q}$ thus $f\notin\mathfrak{q}$ since $\mathfrak{q}$ is a pure ideal. Hence, $\Min(A)\cap D(f)$ is an open in the pure topology for all $f\in A$. The converse implication is deduced from the well known fact that a ring $A$ is a p.p. ring iff $A$ is a reduced mp-ring and $\Min(A)$ is Zariski compact, see \cite[Theorem 4.1]{Tarizadeh p.p. rings} or \cite[Proposition 3.4]{Vasconcelos}. $\Box$ \\

In \cite[Chap. 7, Example 36]{Borceux}, a non-commutative ring is given which has a two-sided maximal ideal whose pure part is not a purely-maximal ideal. But it seems that finding a ``commutative'' ring with a purely-prime ideal which is not a purely-maximal ideal is not easy at all. Hence, this leads us to propose the following conjecture. \\

\begin{conjecture}\label{Conjecture I} Every purely-prime ideal of a commutative ring is purely-maximal. \\
\end{conjecture}

Prove or disprove of the above conjecture would be certainly a nontrivial result. It seems to us that the disproving of the above conjecture looks more likely. But finding a counterexample is not as easy as one may think at first, because we observed that this conjecture holds for both Gelfand rings and reduced mp-rings. Most of the rings which appear in commutative algebra (and algebraic geometry) are either Gelfand rings or reduced mp-rings. \\

\section{Semi-Noetherian rings}

In this section we introduce and study the new notion of semi-Noetherian ring. \\

\begin{definition} We call a ring $A$ a \emph{semi-Noetherian ring} if every pure ideal of $A$ is a finitely generated ideal. \\
\end{definition}

Every Noetherian ring is obviously a semi-Noetherian ring, but the converse is not true. As an example, if $k$ is a field then the polynomial ring $k[x_{1},x_{2},x_{3},...]$ is a semi-Noetherian ring which is not a Noetherian ring. In fact, every domain and every local ring are semi-Noetherian rings. If $R$ is a non-zero ring (i.e. $0\neq1$) then the ring $A=\prod\limits_{i\geq1}R$ is not a semi-Noetherian ring, because $I=\bigoplus\limits_{i\geq1}R$ is a pure ideal of $A$ which is not a finitely generated ideal. If $A$ is a semi-Noetherian ring and $I$ a pure ideal of $A$, then $A/I$ is a semi-Noetherian ring. If $I=(f_{1},...,f_{n})$ is a finitely generated and pure ideal of a ring $A$ then there exists some $g\in I$ such that $f_{i}=f_{i}g$ for all $i$. It follows that $(1-g)I=0$. Hence, $g$ is an idempotent and $I=Ag$. \\

\begin{theorem}$($Cohen type theorem$)$\label{Theorem Cohen type} If every purely-maximal ideal of a ring $A$ is finitely generated, then $A$ is a semi-Noetherian ring.  \\
\end{theorem}

{\bf Proof.} If $P$ is a purely-prime ideal of $A$ then there exists a purely-maximal ideal $M$ of $A$ such that $P\subseteq M$. By the above argument, there exists an idempotent $f\in M$ such that $M=Af$. But $f\in P$ and so $P=M$. Then we prove that every pure ideal of $A$ is a finitely generated ideal. Let $\mathcal{S}$ be the set of all pure ideals of $A$ which are not finitely generated. It suffices to show that $\mathcal{S}=\emptyset$. If not, then by the Zorn's Lemma, it has a maximal element $J$. We show that $J$ is a purely-prime ideal of $A$. Clearly $J\neq A$. Let $I_{1}$ and $I_{2}$ be two pure ideals of $A$ such that $I_{1}I_{2}\subseteq J$. Suppose $I_{1}\nsubseteq J$ and $I_{2}\nsubseteq J$. Note that $J+I_{1}$ and $J+I_{2}$ are pure ideals. It follows that $J+I_{1}=(f_{1},...,f_{n})$ and $J+I_{2}=(g_{1},...,g_{m})$ are finitely generated ideals. If $h\in J$ then there exists some $h'\in J$ such that $h=hh'$. We may write $h=\sum\limits_{i=1}^{n}r_{i}f_{i}$ and $h'=\sum\limits_{k=1}^{m}r'_{k}g_{k}$. It follows that $h=\sum\limits_{i,k}r_{i}r'_{k}f_{i}g_{k}$. This yields that $J=(f_{i}g_{k}: i\in\{1,...,n\}, k\in\{1,...,m\})$ is a finitely generated ideal. But this is a contradiction and we win. $\Box$ \\

\textbf{Acknowledgements.} We would like to thank the referee for very careful reading of the paper and for his/her valuable suggestions and comments which improved the paper. \\

\end{document}